\documentclass{article}
\usepackage{stmaryrd}
\usepackage{amsfonts}
\usepackage{graphicx}
\usepackage{times}
\usepackage{mathptm}
\usepackage{amsmath}

\begin{document}

\title{\ \\ \LARGE\bf Using Gradual Numbers to Analyze Non-Monotonic Functions of Fuzzy Intervals}

\author{\normalsize
         Elizabeth Untiedt$^{1}$\quad
         Weldon Lodwick$^{1}$\\
         \normalsize\sl
         1.Department of Mathematical Sciences, University of Colorado at Denver\\     
        \normalsize\sl
        Email: euntiedt@math.cudenver.edu, wlodwick@math.cudenver.edu
        }

\maketitle

\begin{abstract}
Gradual numbers have been introduced recently as a means of extending standard interval computation methods to fuzzy intervals.  The literature treats monotonic functions of fuzzy intervals.  
In this paper, we combine the concepts of gradual numbers and optimization, which allows for the evaluation of any differentiable function on
 fuzzy intervals, with no monotonicity requirement.    
\end{abstract}


\section{Introduction}
The evaluation of a function on a fuzzy interval has long been accomplished
by evaluating the function on each $\alpha$-cut of the the fuzzy interval.  
Fortin, Dubois, and Fargier recently introduced the 
concept of gradual numbers \cite{fortin07, dubois04}.  A fuzzy interval
can be represented as an interval of gradual numbers, which allows the 
extension of standard interval computation methods to fuzzy intervals.

Previously published results have extended the practice of analyzing monotonic
functions on intervals by evaluating the interval extreme points \cite{fortin07}.  In the 
present paper, we employ optimization techniques on gradual numbers to 
evaluate differentiable functions on fuzzy intervals with no monotonicity
requirements.

\section{Classical Interval Computation}

Classical interval analysis can be viewed as seeking upper and lower bounds for the value of a operator on  
interval-valued arguments.  For example, when the operator is a real-valued 
function,  $f: \textbf{R}^n \mapsto \textbf{R}$ and $x_1 \in [x_1^-, x_1^+]$,..., $x_n \in [x_n^-, x_n^+]$, the goal is to 
find $z^- $ and $z^+$ such that $z$=$f(x_1,...,x_n)$ lies in the 
interval $[z^-,z^+]$ for every $x \in \mathcal{X} = \times_i [x_i^-,x_i^+]$.

Warmus \cite{warmus56}, Sunaga \cite{sunaga58}, and Moore \cite{moore62} 
first proposed what Lodwick 
\cite{lodwick07} calls an axiomatic approach to interval arithmetic, 
which defines the following rules for basic interval arithmetic:

For $x \in [x^-,x^+]$ and $y \in [y^-,y^+]$: 
\begin{enumerate}
\item $x + y = [x^- + y^-, x^+ + y^+]$,
\item $x - y = [x^- - y^+, x^+ - y^-]$,
\item $x \times y = [\min\{x^- \times y^-, x^- \times y^+, x^+ \times y^-, x^+ \times y^+ \},\max\{x^- \times y^-, x^- \times y^+, x^+ \times y^-, x^+ \times y^+ \}]$,
\item $x \div y = [\min\{x^- \div y^-, x^- \div y^+, x^+ \div y^-, x^+ \div y^+ \},\max\{x^- \div y^-, x^- \div y^+, x^+ \div y^-, x^+ \div y^+ \}]$, $0 \notin [y^-,y^+]$. 
\end{enumerate}

Axiomatic interval arithmetic assumes that all variables are independent, 
including duplicate copies of the same variable.  This can lead to 
overestimations of the bounds of functions.  
 Consider, for example,
\[f(x) = x - x,~x \in [0,1].\]
 
Using the axiomatic approach yields the result 
\[ [1,0]-[1,0] = [0-1,1-0] = [-1,1].\]

If, however, we simplify the function {\em before} applying axiomatic 
interval analysis, we get $f(x)= x-x=0$, which indicates that $[z^-,z^+] = 
[-1,1]$ is a 
gross over-estimation.

  Simplifying functions to reduce multiple copies of a variable can eliminate the overestimation, but is not always possible.  For example, there is no way to express $f(x) = xe^x$ with only one copy of the variable $x$.

When the function in question is monotonic, we can appeal to a well-known proposition 
to find exact bounds for $[z^-,z^+]$ without overestimation.  Before stating 
this result, we provide some preliminary definitions.  A 
{\it configuration}, $\omega$, is an $n$-tuple which lies in $\mathcal{X}=\times_i[x_i^-,x_i^+]$.  The set of {\it extreme configurations}, $\mathcal{H}$, are $n$-tuples
 in which every element takes on the value of an interval endpoint.  Of course, there are $2^n$ elements in $\mathcal{H}$.    Now 
consider the following proposition, which indicates what is known 
as the vertex method of interval analysis in the fuzzy set community.

\textit{Proposition 1:}  Let $ x = (x_1,...,x_n)$ be a tuple of $n$ real interval-valued variables such that $x_i \in [x_i^-,x_i^+]$, and let $z=f(x_1,...,x_n)=[z^-,z^+].$  
If $f$ is continuous  and locally monotonic with respect to each argument, then \[z^- =  \displaystyle{\min_{\omega \in \mathcal{H}}}f(\omega)\]
and \[z^+ =  \displaystyle{\max_{\omega \in \mathcal{H}}}f(\omega).\]

This follows from the fact that $\mathcal{X}=\times_i [x_i^-,x_i^+]$ is a closed subset 
of \textbf{$R^n$}, and from the Weierstrass theorem that a continuous 
function attains its maximum and minimum on a compact set.
  The vertex method, while effective, is limited in scope.  
We are still left with over-estimations or underestimations on non-monotonic 
functions.  Take, for example,     
                                                                                
\begin{equation}\label{eq:xminus1}
f(x) = x(1-x)~~for~~x \in [0,1].\end{equation}                                             
If we consider the two instantiations of $x$ to be independent, 
$\mathcal{H}$ contains the four
 elements $\{ (0,0), (0,1), (1,0), (1,1)\}$, and the result is 
\[z^- = f(\omega)^- = \min_{\omega \in \mathcal{H}}f(\omega) = f(1,1) = f(0,0) = f(0,1) = 0,\]
and 
\[z^+ = f(\omega)^+ = \max_{\omega \in \mathcal{H}}f(\omega) = f(1,0) = 1.\]

Inspection reveals that there is no single value of $x$ in the interval 
$[0,1]$ for which $f(x) = x(x-1)$ has a value of 1, so 
$[z^-,z^+] = [0,1]$  is an unsatisfactory solution.  We might, on the 
other hand, consider the two instantiations 
of $x$ to be the same variable.  Then $\mathcal{H}$ contains 
only the two elements $\{(0,0), (1,1)\}$, and the result is      
\[z^- = f(\omega)^- = \min_{\omega \in \mathcal{H}}f(\omega) = f(1,1) = f(0,0) = 0,\]            
and                                                                             \[z^+ = f(\omega)^+ = \max_{\omega \in \mathcal{H}}f(\omega) = f(1,1) = f(0,0) = 0.\]

But $[z^-,z^+] = [0,0]$ is also wrong because zero is not the only possible value of $f$.  Consider, for example, $f(\frac{1}{4}) = \frac{3}{16}.$  We require a broader tool than the vertex method for evaluation of
 non-monotonic functions of interval-valued variables with dependencies.

\section{Optimization and Interval Analysis}

Observe that the values we seek in interval analysis,  $z^- = \min_{x \in \mathcal{X}} f(x)$ and  $z^+ = \max_{x \in \mathcal{X}} f(x)$ are the answers to the optimization problems: 
\begin{eqnarray}
\mbox{Minimize}~&f(x)&\label{eq:MP1obj}\\
\mbox{subject~to}~& x_i &\ge x_i^-,~~~\mbox{for}~i=1,...,n\label{eq:MP1con}\\
                 & -x_i &\ge -x_i^+,~~~\mbox{for}~i=1,...,n,\nonumber
\end{eqnarray} 

and
 \begin{eqnarray}
\mbox{Maximize}~&f(x)&\label{eq:MP2obj}\\
\mbox{subject~to}~& x_i &\ge x_i^-,~~~\mbox{for}~i=1,...,n\label{eq:MP2con}\\
                 & -x_i &\ge -x_i^+,~~~\mbox{for}~i=1,...,n.\nonumber
\end{eqnarray} 

Lodwick proposed this idea in the nineties with his constrained interval 
analysis (\cite{lodwick99}), and more recently Kreinovich (\cite{kreinovich04}) 
has applied it to technical problems with partial information.

Viewing the interval analysis problem as a mathematical programming problem puts at our disposal the rich optimization theory of the last century, which leads to the following proposition.

\textit{Proposition 2:}  Let $ x = (x_1,...,x_n)$ be a tuple of $n$ 
variables such that $x_i \in [x_i^-,x_i^+]$ 
and let $z=f(x_1,...,x_n)$ such that  $f$ is continuously differentiable with respect to each $x_i$. 
Now consider the set of KKT configurations, $\mathcal{K}^-_{f}$ and 
$\mathcal{K}^+_{f}$, to be configurations which satisfy the Karush-Kuhn-Tucker 
conditions for (\ref{eq:MP1obj})-(\ref{eq:MP1con}) and (\ref{eq:MP2obj})-
(\ref{eq:MP2con}) respectively.  Then 
\[z^- =  \min_{\omega \in \mathcal{K^-}}f(\omega)\] 
and \[z^+ =  \max_{\omega \in \mathcal{K}+}f(\omega).\]
                         
\textit{Proof:}      Finding the minimum function value for $f(x_1,...,x_n)$ where $x \in \mathcal{X}=\times_i[x_i^-,x_i^+]$ can be stated as a mathematical programming problem:
\begin{eqnarray}
\mbox{Minimize}~z&:&f(x)\label{eq:MP3obj}\\
\mbox{subject~to}~c(i)&:&x_i-x_i^- \ge 0,~~~\mbox{for}~i=1,...,n\label{eq:MP1thrun}\\
                 c(n+i)&:&x_i^+-x_i\ge 0,~~~\mbox{for}~i=1,...,n,\label{eq:MPnplus1}
\end{eqnarray} 
 where $x$ is a $1 \times n$ vector.

Since constraints $c(i)$ and $c(n+i)$ are lower and upper bounds for $x_i$, 
at most one at a time of $c(i)$ and $c(n+i)$ will be active.  
We observe that our mathematical programming problem satisfies the linear independence constraint qualifications (i.e. the gradients of the active constraints are linearly independent),
which means that the KKT conditions are necessary conditions for optimality.  The details of 
Karush-Kuhn-Tucker conditions are found in any standard optimization textbook, for example \cite{nocedalandwright}.  

Now any optimal point will satisfy the following:
\begin{eqnarray}
\bigtriangledown_x f(x)- \lambda_1(x_i-x_i^-) &-&... \\ 
... -\lambda_n(x_n-x_n^-) + \lambda_{n+1}(x_1-x_1^+) &+& ... \\ 
... +\lambda_{2n}(x_b-x_n^+) &=&0 \label{eq:KKT1}\\
\lambda_i\ge 0,~~~&&~i=1,...,2n\label{eq:KKT2}\\
x_i-x_i^- \ge 0,~~~&&~i=1,...,n\label{eq:KKT3}\\
-x_i + x_i^+ \ge 0,~~~&&~i=1,...,n\label{eq:KKT4}\\
\lambda_i(x_i-x_i^-) = 0,~~~&&~i=1,...,n\label{eq:KKT5}\\
\lambda_n+i(-x_i+x_i^+) = 0,~~~&&~i=1,...,n
\label{eq:KKT6}.
\end{eqnarray}

The set of points satisfying (\ref{eq:KKT1})-(\ref{eq:KKT5}) are exactly the points in $\mathcal{K}^-_f$.  Note that $\mathcal{X}= \times_i [x_i^-,x_i^+]$ is a 
closed subset of \textbf{$R^n$}, so is a compact set.  The Weierstrass theorem guarantees that a continuous function attains a maximum and minimum on a compact set, so we know $f(x)$ has a minimum on $\mathcal{X}$,
which has to be $\min_{x \in \mathcal{K}^-_f} f(x).$  
  
To find the maximum, we replace $f(x)$ in (\ref{eq:MP3obj}) with $-f(x)$, and find the KKT points, which are exactly $\mathcal{K}^+_f$.

\subsection{Examples of Interval Analysis via Optimization}
For the sake of illustration, recall example (\ref{eq:xminus1}):
\begin{equation}\label{eq:x-xsquared}
f(x) = x(1-x),~x \in [0,1].\end{equation}

The Karush-Kuhn-Tucker conditions for (\ref{eq:x-xsquared}) are
\begin{eqnarray}
-2x+1-\lambda_1+\lambda_2 &=&0 \label{eq:KKT1ex}\\
\lambda_i &\ge& 0,~j=1,2\label{eq:KKT2ex}\\
x &\ge& 0,\label{eq:KKT3ex}\\
 1-x &\ge& 0,\label{eq:KKT4ex}\\
\lambda_1 x&=&0\label{eq:KKT5ex}\\
\lambda_2 (1-x) &=&0\label{eq:KKT6ex}.
\end{eqnarray}

$\mathcal{K}^-_f$ consists of the three points which satisfy (\ref{eq:KKT1ex})-(\ref{eq:KKT6}): 
\begin{itemize} 
\item $x = 0; ~\lambda = (1,0);~ f(x) =0$
\item  $x = 1; ~\lambda = (0,1);~ f(x) =0$
\item  $x = \frac{1}{2}; ~\lambda = (0,0);~ f(x) =\frac{1}{4}.$
\end{itemize}

 Now  $z^- = \min_{x \in \mathcal{K}_f^-} f(x)$ = 0.  In this example, it happens that $\{x \in \mathcal{K}_f^+\}=\{x \in \mathcal{K}_f^-$\}, so
$z^+ = \max_{x \in \mathcal{K}_f^+} f(x) = \frac{1}{4}$.  

Optimization yields the desired result, $[z^-,z^+]=[0,\frac{1}{4}]$ without overestimation or underestimation.  
It should be noted that this accuracy is not without its cost.  The vertex method tells us that, 
for monotonic functions $f:\textbf{R}^n \mapsto \textbf{R}$, function extrema will always occur at interval endpoints, so 
$|\mathcal{H}| = 2^n$.  In contrast, the number of KKT configurations of an interval constrained mathematical programming problems
 will sometimes be fewer than the number of extreme configurations, and will occasionally be many, many more.  In 
the worst case, the number of points satisfying the KKT conditions will be infinite.

\section{Gradual Numbers and Fuzzy Interval Analysis}

Up to this point, we have considered intervals of real numbers, which provide a Boolean model for uncertainty.  
Numbers inside the interval are entirely possible, while numbers outside the interval are impossible.  It is often the case, 
however, that the sets we wish to model are not so crisply defined.  The
need frequently arises to move gradually from impossible to possible,
with values on the outskirts of the interval being somewhat possible.  
Fuzzy set theory accomplishes this by making the bounds of the interval softer, 
so that the transition from an impossible value to an entirely possible value is gradual.  

A fuzzy interval, $\tilde{W}$ can be defined by a membership function, 
$\mu_{\tilde{W}}: \textbf{R} \mapsto (0,1]$, where  $\mu_{\tilde{W}}(x)$ 
indicates to what extent $x$ is a member of $\tilde{W}$.  A standard fuzzy interval is the $L-R$ fuzzy interval, which is defined by a
membership function $\mu_{\tilde{W}}$ and reference functions $L$ and $R$.  $L:[0,\inf)\rightarrow[0,1]$ and $R:[0,\inf)\rightarrow[0,1]$ are upper semi-continuous and strictly decreasing
in the range (0,1], and $\mu_{\tilde{W}}$ is defined as follows:
\begin{eqnarray}                                                              
\mu_{\tilde{W}}(x) &=& \begin{cases}                                                  
1&\mbox{for~}x\in[w^-,w^+],\\
L(\frac{w^--x}{\alpha_W})&\mbox{for~} x<w^-,\\       
R(\frac{x-w^+}{\beta_W})&\mbox{for~}x>w^+.             
\end{cases}                   
\end{eqnarray}

Fortin, Dubois, and Fargier \cite{fortin07} recently proposed that a fuzzy interval can be represented as a set of gradual numbers that lie between two gradual number endpoints 
in the same way that a real interval is a set of real numbers that lie between two real endpoints.  

{\textit Definition 1:}  A \textit{gradual real number} (or gradual number), $\tilde{r}$, is defined by an assignment function 
$A_{\tilde{r}}:(0,1] \mapsto {\textbf{R}}$.  It can be understood as a real value parametrized  by $\alpha \in (0,1]$.  Then for each 
$\alpha$, a real value $r_{\alpha}$ is given by $A_{\tilde{r}}(\alpha)$.  Let 
$G$ be the set of all gradual real numbers.  

We can describe a fuzzy interval $\tilde{W}$ by an ordered pair of two gradual numbers, $[\tilde{w}^-,\tilde{w}^+]$.  Since a fuzzy interval is a fuzzy set with an upper-semicontinuous
membership function, we require that $A_{\tilde{w}^-}$ be increasing, that  $A_{\tilde{w}^+}$ be decreasing, and that  $A_{\tilde{w}^-}(1) \le  A_{\tilde{w}^+}(1)$.  
For our results, we will additionally assume that  $A_{\tilde{w}^-}$ and  $A_{\tilde{w}^+}$ are continuous and that  $A_{\tilde{w}^-}(0)$ and  $A_{\tilde{w}^+}(0)$ are defined.      

Before considering arithmetic operations on {\em intervals} of gradual numbers, we need a definition of arithmetic operations on gradual numbers themselves.  If $\tilde{r}$ 
and $\tilde{s}$ are two gradual numbers, then $\tilde{r} + \tilde{s}$ is 
defined by its membership function $\mathcal{A}_{\tilde{r} + \tilde{s}} = 
\mathcal{A}_{\tilde{r}} + \mathcal{A}_{\tilde{s}}$.  Subtraction, multiplication, division, minimum and maximum are similarly defined.  In addition, 
$\tilde{r}$ is said to be less than $\tilde{s}$ if $\mathcal{A}_{\tilde{r}}(\alpha) < \mathcal{A}_{\tilde{s}}(\alpha)$ for all $\alpha \in [0,1]$. Greater than and equality are similarly defined.  We would like to note that the set of gradual numbers is not a fully ordered set.     

Beyond arithmetic operations on gradual numbers, we will find it useful to
 extend a real-valued function $f(x_1,...,x_n): R^n \mapsto R$ 
to a gradual function, $\dot{f}(\tilde{x}_1,...,\tilde{x}_n): G^n \mapsto G$. 
We define the assignment function of the value of $\dot{f}$ as follows:  
\begin{eqnarray}\label{eq:gradualextension}
\mathcal{A}_{\dot{f}(\tilde{r}_1,...,\tilde{r}_n)}(\alpha)&=& 
f(\tilde{r}_1(\alpha)...,\tilde{r}_n(\alpha)),\\
& & \forall \alpha  \in  [0,1]                            
\end{eqnarray}

A fuzzy interval $\tilde{W} = [\tilde{w}^-,\tilde{w}^+]$ is the set of all gradual numbers, $\tilde{w}$, for which $\tilde{w}^- \le \tilde{w} \le \tilde{w}^+$.  If $\tilde{W}$ is an L-R fuzzy interval, then $\tilde{w}^-$ and $\tilde{w}^+$
are defined by the following assignment functions:
\begin{eqnarray}
\mathcal{A}_{\tilde{W}^-}(\alpha) &=& w^- - L^{-1}(\alpha) \alpha_W\\
\mathcal{A}_{\tilde{W}^+}(\alpha) &=& w^+ + R^{-1}(\alpha) \beta_W,
\end{eqnarray}

where $L^{-1}$ (respectively $R^{-1}$) is the inverse of $L$ 
(respectively $R$) in the part of its domain where it is positive.

\subsection{Earlier Results on Applying Interval Analysis to Fuzzy Intervals}

The representation of a fuzzy interval as a crisp interval of gradual numbers 
allows us to use tools previously developed for real interval analysis in the
analysis of fuzzy intervals.

Fortin, Dubois and Fargier extend the vertex method of interval analysis to intervals of gradual numbers \cite{fortin07}.  They define
a fuzzy configuration, $\Omega$, which is a tuple of gradual numbers, $(\tilde{r}_1,...,\tilde{r}_n)$ as follows:

Let fuzzy intervals $\tilde{x_1}$...$\tilde{x}_n$ be described by their gradual endpoints, 
$ [\tilde{x}_1^-,\tilde{x}_1^+]...[\tilde{x}_n^-,\tilde{x}_n^+]$.  A {\it configuration}, $\Omega$, is an $n$-tuple of gradual numbers which lies in 
$\mathcal{\tilde{X}}=\times_i[\tilde{x}_i^-,\tilde{x}_i^+]$.  The set of 
{\it extreme fuzzy configurations}, $\tilde{\mathcal{H}}$, are $n$-tuples in 
which every 
element takes on the value of an interval endpoint.   A fuzzy configuration 
or extreme configuration evaluated at $\alpha$ 
is defined as the tuple of each of its gradual components evaluated at 
$\alpha$.  That is $\Omega(\alpha) = (\tilde{r}_1(\alpha), ... , 
\tilde{r}_n(\alpha)).$  We now turn to the vertex method for fuzzy interval analysis.                                                                                                       
                                                                                                                                               
\textit{Proposition 3:} Let $\tilde{x} = (\tilde{x}_1,...,\tilde{x}_n)$ be a 
tuple of $n$ fuzzy interval-valued variables such that 
$\tilde{x}_i \in [\tilde{x}_i^-,
\tilde{x}_i^+]$.  Let $\dot{f}: G^n \mapsto G$ be the gradual extension of 
$f: {\textbf R^n} \mapsto {\textbf R}$  as defined in 
(\ref{eq:gradualextension}),
and define the interval $[\tilde{z}^-,\tilde{z}^+]$ to be the smallest 
interval of gradual numbers that contains 
$\tilde{z}=\{\dot{f}(\tilde{x}_1,...,\tilde{x}_n)\}$ for every $\tilde{x} \in \tilde{\mathcal{X}}=\times_i [\tilde{x}_i^-,\tilde{x}_i^+]$.                                                                                                                                           
If $f$ is continuous and locally monotonic with respect to each argument, then                                                                
\[\tilde{z}^- =  \displaystyle{\dot{\min}_{\Omega \in \tilde{\mathcal{H}}}}\dot{f}(\Omega),\]                                                                               
and \[\tilde{z}^+ =  \displaystyle{\dot{\max}_{\Omega \in \tilde{\mathcal{H}}}}\dot{f}(\Omega).\]                                                                          
                                                                                                                                               
The reader is referred to \cite{fortin07} for the proof. 

\begin{figure}[htp]
\centerline{\includegraphics[width=3.0in]{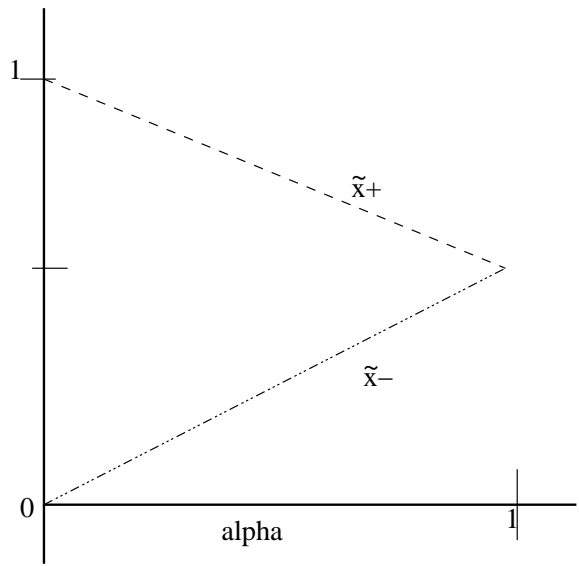}}
\caption{$\tilde{x}=[\tilde{x}^-,\tilde{x}^+] = [\frac{1}{2}\alpha,
1-\frac{1}{2}\alpha]$.}\label{fig:fuzzyx}
\end{figure}

An example of the vertex method for intervals of gradual numbers will be instructive.
Let $\tilde{x} = [\tilde{x}^-,\tilde{x}^+]$ 
be a fuzzy interval defined by the following gradual endpoints:
\begin{eqnarray}\label{eq:assignment1}
\tilde{x}^-(\alpha)& =& \frac{1}{2}\alpha, ~~\alpha \in [0,1] \\
\tilde{x}^+(\alpha)& =& 1-\frac{1}{2}\alpha, ~~\alpha \in [0,1], \\   
\end{eqnarray}
as in Figure \ref{fig:fuzzyx}.  

Let $\dot{f}=\tilde{x}^2$ for $\tilde{x} \in
\tilde{X}$.  To apply the vertex method, we determine that there are two elements 
in $\tilde{\mathcal{H}}$: $\tilde{x}^-(\alpha) = \frac{1}{2}\alpha, ~~\alpha \in [0,1]$ 
and $\tilde{x}^+(\alpha) = 1-\frac{1}{2}\alpha, ~~\alpha \in [0,1]$. Now                   
\begin{eqnarray}                                                                         
\dot{f}(\tilde{x}^-(\alpha))& =& \frac{1}{4}\alpha^2,~\alpha \in [0,1] \\          
\dot{f}(\tilde{x}^+(\alpha))&=& (1-\frac{1}{2}\alpha)^2 \\
              &=& 1- \alpha + \frac{1}{4}\alpha^2,~\alpha \in [0,1].      
\end{eqnarray}                                                          

\begin{figure}                                                                 \centerline{\includegraphics[width=3.0in]{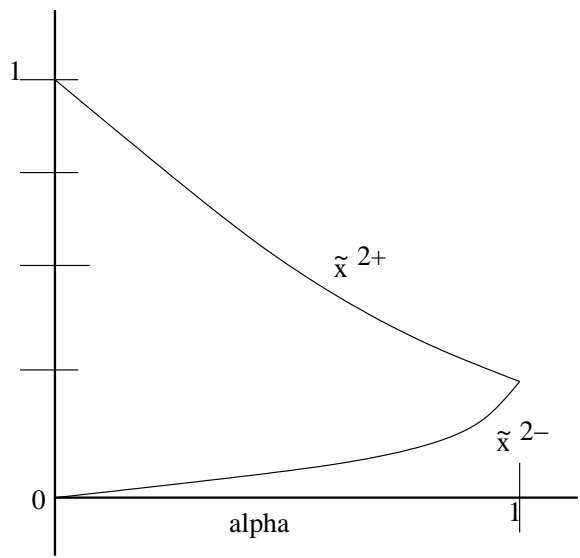}}
\caption{$\tilde{x}^2 = [\tilde{z}^-,\tilde{z}^+]=[\frac{1}{4}\alpha^2,
1-\alpha+\frac{1}{4}\alpha^2]$.}\label{fig:fuzzyx2}
\end{figure}

It is clear that $\dot{f}(\tilde{x}^-(\alpha)) \le \dot{f}(\tilde{x}^+(\alpha))$ for all
 $\alpha \in [0,1]$, so $z^- = \dot{f}(\tilde{x}^-)$ and $z^+ = \dot{f}(\tilde{x}^+)$.  The fuzzy
interval $\tilde{z} = [\tilde{z}^-,\tilde{z}^+]$ is shown in Figure \ref{fig:fuzzyx2}.
 
\begin{figure}
\centerline{\includegraphics[width=3.0in]{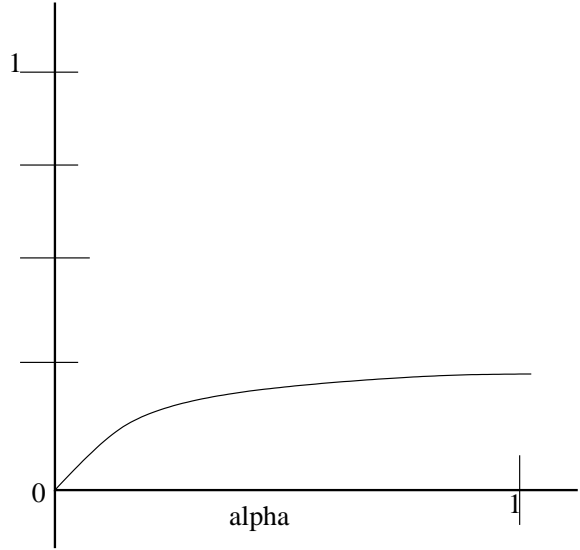}}
\caption{$\tilde{z}^-=\tilde{z}^+ = \frac{1}{2}\alpha-\frac{1}{4}(\alpha)^2$.}
\label{fig:fuzzyxminusx2}
\end{figure}

The vertex method, while effective, is limited in scope.  We are still left              
with over-estimations, underestimations, or meaningless results 
on non-monotonic                               
functions.  For example, let $\tilde{X}$ be defined as above, and
let $\dot{f}$ be the gradual extension of $f(x)=x(1-x)$ for $\tilde{x} \in \tilde{X}$.  To 
apply the vertex method, we determine that there are two elements in 
$\tilde{\mathcal{H}}$: $\tilde{x}^-(\alpha) = \frac{1}{2}\alpha, ~~\alpha \in [0,1]$ and
$\tilde{x}^+(\alpha)= 1-\frac{1}{2}\alpha, ~~\alpha \in [0,1]$. Now
\begin{eqnarray}
\dot{f}(\tilde{x}^-(\alpha))& =& \frac{1}{2}\alpha (1-\frac{1}{2}\alpha),~\alpha \in [0,1] \\
              &=& \frac{1}{2}\alpha - \frac{1}{4}\alpha^2,~\alpha \in [0,1] \\ 
              &=& \dot{f}(\tilde{x}^+(\alpha)).
\end{eqnarray}
Then $\tilde{z}^-$=$\tilde{z}^+$, so the fuzzy interval is not an interval at all, but a single gradual number defined by the assignment function 
$\mathcal{A}_{\tilde{z}}(\alpha)=\frac{1}{2}\alpha-\frac{1}{4}(\alpha)^2$, as in 
Figure \ref{fig:fuzzyxminusx2}.  
 
\section{Main Results} 
Fortin, Dubois, and Fargier have used the concept of gradual numbers to extend
interval techniques to the evaluation of monotonic functions of fuzzy intervals.
In this section we expand the use of gradual numbers.  We employ
optimization theory in conjunction with gradual numbers to evaluate
continuously differentiable non-monotonic functions of fuzzy-interval valued 
variables.  
 
\subsection{Optimization on Fuzzy Intervals}

Let $\dot{f}: G^n \mapsto G$ be the gradual extension of a continuously 
differentiable function 
$f: {\textbf R^n} \mapsto {\textbf R}$ as defined in 
(\ref{eq:gradualextension}),
and define the interval $[\tilde{z}^-,\tilde{z}^+]$ to be the smallest interval of gradual number that contains 
$\tilde{z}=\{\dot{f}(\tilde{x}_1,...,\tilde{x}_n)\}$ for every 
$\tilde{x} \in \tilde{\mathcal{X}}=\times_i [\tilde{x}_i^-,\tilde{x}_i^+]$.                                                           \

Observe that the values we seek,  
$\tilde{z}^- = \dot{\min}_{\tilde{x} \in \tilde{\Omega}} \dot{f}(x)$ and  
$\tilde{z}^+ = \dot{\max}_{\tilde{x} \in \tilde{\Omega}} \dot{f}(\tilde{x})$, are 
the answers to the optimization problems                                        \begin{eqnarray}                                                                         
\mbox{Minimize}~&\dot{f}(\tilde{x})&\label{eq:FP1obj}\\                                                
\mbox{subject~to}~& \tilde{x}_i &\ge \tilde{x}_i^-,~~~\mbox{for}~i=1,...,n\label{eq:FP1con}\\            
                 & \tilde{x}_i &\le \tilde{x}_i^+,~~~\mbox{for}~i=1,...,n,\nonumber                          
\end{eqnarray}                                                                           
and
\begin{eqnarray}                                                                                                                                             
\mbox{Maximize}~&\dot{f}(\tilde{x})&\label{eq:FP2obj}\\                                                                                                                 
\mbox{subject~to}~& \tilde{x}_i &\ge \tilde{x}_i^-,~~~\mbox{for}~i=1,...,n
\label{eq:FP2con}\\                                                                                    
            & \tilde{x}_i &\le \tilde{x}_i^+,~~~\mbox{for}~i=1,...,n,\nonumber 
\end{eqnarray}   

We now define the set of KKT configurations, $\tilde{\mathcal{K}}^-_{f}$ and     
$\tilde{\mathcal{K}}^+_{f}$,                                                    
to be configurations which satisfy the Karush-Kuhn-Tucker conditions for        
(\ref{eq:FP1obj})-(\ref{eq:FP1con}) and (\ref{eq:FP2obj})-(\ref{eq:FP2con}) 
respectively. 

Explicitly, $\tilde{\mathcal{K}}^+_{f}$ is the set of every configuration $x \in \tilde{\mathcal{X}}$ 
which satisfies the following:
\begin{eqnarray}                                                                                                                                             
\dot{g}(\tilde{x})- \lambda_1(\tilde{x}_i-\tilde{x}_i^-) -  ...&& \\
...-\lambda_n(\tilde{x}_n-\tilde{x}_n^-) + \lambda_{n+1}(\tilde{x}_1-\tilde{x}_1^+)  + ...&&\\
... +\lambda_{2n}(\tilde{x}_i-\tilde{x}_n^+) &=&0 \label{eq:FKKT1}\\                        
\lambda_i\ge 0,&&~i=1,...,2n\label{eq:FKKT2}\\                                                                                                     
\tilde{x}_i-\tilde{x}_i^- \ge 0,&&~i=1,...,n\label{eq:FKKT3}\\      
-\tilde{x}_i + \tilde{x}_i^+ \ge 0,&&~i=1,...,n\label{eq:FKKT4}\\   \lambda_i(\tilde{x}_i-\tilde{x}_i^-) = 0,&&~i=1,...,n
\label{eq:FKKT5}\\                                                                             
\lambda_n+i(-\tilde{x}_i+\tilde{x}_i^+) = 0,&&~i=1,...,n      
\label{eq:FKKT6},                                                              
\end{eqnarray}     
where $\dot{g}$ is the gradual extension of $\bigtriangledown
f(x_1,...,x_n)$.  We define $\tilde{\mathcal{K}}^+(\alpha) = \{\tilde{x}(\alpha)|
\tilde{x} \in \tilde{\mathcal{K}}^+$.  $\tilde{\mathcal{K}}^-(\alpha)$ is similarly defined.

\textit{Proposition 4:}~~Let $ \tilde{x} = (\tilde{x}_1,...,\tilde{x}_n)$ be a tuple of $n$ fuzzy 
interval-valued variables such that 
$\tilde{x}_i$ is any gradual number in the interval defined by gradual numbers 
$[\tilde{x}_i^-,\tilde{x}_i^+]$.  Let $[\tilde{z}^-,\tilde{z}^+]$ be the
smallest fuzzy interval that contains $\tilde{z} = \dot{f}(\tilde{x}_1,...,\tilde{x}_n)$ for every $\tilde{x} \in \tilde{\mathcal{X}}$,
where  $\dot{f}$ is the extension of a continuously 
differentiable real-valued function $f:R^n \mapsto R$.  Then 
\[\tilde{z}^- =  \dot{\min}_{\tilde{\omega} \in \tilde{\mathcal{K^-}}}\dot{f}(\tilde{\omega})\]                                                    and
\[\tilde{z}^+ =  \dot{\max}_{\tilde{\omega} \in \tilde{\mathcal{K^-}}}\dot{f}
(\tilde{\omega}).\]  

\textit{Proof:}  Finding the minimum function value for $\dot{f}(\tilde{x}_1,...,\tilde{x}_n)$ where 
$\tilde{x} \in \tilde{\mathcal{X}}=\times_i[\tilde{x}_i^-,\tilde{x}_i^+]$ can be stated as a mathematical programming problem:       
\begin{eqnarray}                                                                         
\mbox{Minimize}~\tilde{z}&:&\dot{f}(\tilde{x})\label{eq:FP3obj}\\                                              
\mbox{subject~to}~c(i)&:&\tilde{x}_i-\tilde{x}_i^- \ge 0,~~~\mbox{for}~i=1,...,n\label{eq:FP1thrun}\\    
                 c(n+i)&:&\tilde{x}_i^+-\tilde{x}_i\ge 0,~~~\mbox{for}~i=1,...,n,\label{eq:FPnplus1}     
\end{eqnarray}                                                                           
 where $\tilde{x}$ is a $1 \times n$ vector of fuzzy interval-valued variables. 

Fix $\alpha \in [0,1]$.  Since constraints $c(i)$ and $c(n+i)$ are lower and upper bounds for 
$\tilde{x}_i$ respectively, at most one at a time will be active.  We observe that our mathematical programming problem satisfies the linear independence 
constraint qualifications (the gradients of the active constraints are linearly 
independent,  which means that the KKT conditions are necessary conditions for optimality.  
Let $g(x)=\bigtriangledown_x f(x)$, and let $\dot{g}(\tilde{x})$ be the extension of 
$g(x)$.

Then for a fixed $\alpha \in [0,1]$, any optimal configuration, $\tilde{x}(\alpha)$ will 
satisfy the following:                                   
\begin{eqnarray}                                                                         
\dot{g}(\tilde{x}(\alpha)- \lambda_1(\tilde{x}_i(\alpha)-\tilde{x}_i^-(\alpha)) -  ...&& \\
... -\lambda_n(\tilde{x}_n(\alpha)-\tilde{x}_n^-(\alpha)) &&\\
+ \lambda_{n+1}(\tilde{x}_1(\alpha)-\tilde{x}_1^+(\alpha))  + ... && \\
... +\lambda_{2n}(\tilde{x}_i(\alpha)-\tilde{x}_n^+(\alpha)) &=& 0 \label{eq:F2KKT1}\\                       
\lambda_i\ge 0,~~~&\mbox{for}&~i=1,...,2n\label{eq:F2KKT2}\\                              
\tilde{x}_i(\alpha)-\tilde{x}_i^-(\alpha) \ge 0,&&~i=1,...,n\label{eq:F2KKT3}\\                              
-\tilde{x}_i(\alpha) + \tilde{x}_i^+(\alpha) \ge 0,&&~i=1,...,n\label{eq:F2KKT4}\\                           
\lambda_i(\tilde{x}_i(\alpha)-\tilde{x}_i^-(\alpha)) = 0,&&~i=1,...,n\label{eq:F2KKT5}\\                     
\lambda_{n+i}(-\tilde{x}_i(\alpha)+\tilde{x}_i^+(\alpha)) = 0,&&~i=1,...,n                                      
\label{eq:F2KKT6}.                                                                        
\end{eqnarray}                                                                           

Note that since $\mathcal{X}(\alpha)= \times_i [\tilde{x}_i^-(\alpha),\tilde{x}_i^+(\alpha)]$ is a                                  
closed subset of \textbf{$R^n$}, it is a compact set.  Weierstrass' theorem guarantees that a continuous function attains a 
maximum and minimum on a compact set, so we know $f(x)$ has a minimum on $\mathcal{X}(\alpha)$.   
which has to be $\dot{\min}_{\tilde{x} \in \tilde{\mathcal{K}}^-_f} \dot{f}(\tilde{x}(\alpha)).$

The set of points satisfying (\ref{eq:F2KKT1})-(\ref{eq:F2KKT6}) for a fixed $\alpha \in [0,1]$ are exactly the points in $\tilde{\mathcal{K}}^-_f(\alpha)$.  The minimum of this set will be  $\dot{\min}_{\tilde{x} \in \tilde{\mathcal{K}}^-_f} \dot{f}(\tilde{x}(\alpha)).$  Taking this result for every $\alpha \in [0,1]$
allows us to  define $\tilde{z}^-$ by its assignment function $\mathcal{A}_{\tilde{z}^-}(\alpha) =  \dot{\min}_{\tilde{x} \in \tilde{\mathcal{K}}^-_f} \dot{f}(\tilde{x}(\alpha)),$ which is the same as $\tilde{z}^- = \dot{\min}_{\tilde{x} \in \tilde{\mathcal{K}}^-_f} \dot{f}(\tilde{x})$.               
                                                                                         
The proof for $\tilde{z^+}$, is anologous, replacing $\dot{g}(x)$ in (\ref{eq:FP3obj}) 
with $-\dot{g}(x)$.

\subsection{Example of Fuzzy Interval Analysis via Optimization}                               
For the sake of illustration, recall our earlier example:                                
                                                                                         
\begin{equation}\label{eq:x-xsquared2}                                                    
f(x) = {x}(1-{x}). \end{equation}

Now let $\dot{f}(\tilde{x})=\tilde{x}(1-\tilde{x})$ be the
gradual extension of $f(x)$, and let $\tilde{x}$ be restricted by the
fuzzy interval $\tilde{X} = [\tilde{x}^-,\tilde{x}^+]$, defined by the 
following gradual endpoints:
\begin{eqnarray}\label{eq:assignment2}
\tilde{x}^-(\alpha)& =& \frac{1}{2}\alpha, ~~\alpha \in [0,1] \\
\tilde{x}^+(\alpha)& =& 1-\frac{1}{2}\alpha, ~~\alpha \in [0,1]. 
\end{eqnarray}

Finding the right endpoint of $\dot{f}(\tilde{x})$ is equivalent to
solving the following optimization problem:
\begin{eqnarray}\label{eq:gradopt}
\mbox{Minimize}~&&\tilde{x}-\tilde{x}^2\label{eq:gradoptobj}\\                                              
\mbox{subject~to}~&&\tilde{x}-\frac{1}{2}\alpha \ge 0 \\  
    &&\frac{1}{2}\alpha -1 - \tilde{x} \ge 0 
\end{eqnarray}

Now $f(x)$ is continuously differentiable, so we can use Proposition 4 to 
evaluate $\dot{f}(\tilde{x})$ on $\tilde{X}$.

The Karush-Kuhn-Tucker conditions for (\ref{eq:x-xsquared}) are
\begin{eqnarray}\label{eq:gradkkt}                                                                         
1-2\tilde{x}-\lambda_1+\lambda_2 &=&0 \label{eq:gradKKT1}\\                                        
\lambda_i &\ge& 0,~j=1,2\label{eq:gradKKT2}\\                                   
\tilde{x}-\frac{1}{2}\alpha &\ge& 0 \label{eq:gradKKT3}\\
1-\frac{1}{2}\alpha - \tilde{x} &\ge& 0\label{eq:gradKKT4} \\
\lambda_1(\tilde{x}-\frac{1}{2}\alpha) &=& 0 \label{eq:gradKKT5}\\
\lambda_2(1-\frac{1}{2}\alpha - \tilde{x}) &=& 0\label{eq:gradKKT6} \\
\lambda_i &\ge& 0, \lambda = 1,2\label{eq:gradKKT7}
 \end{eqnarray}                                                                           
                                                                                         
Now $\tilde{\mathcal{K}}^-_{\dot{f}}$ consists of the three triples, $(\tilde{x}, \lambda_1,\lambda_2)$, which satisfy (\ref{eq:gradKKT1})-(\ref{eq:gradKKT7}):                                                                                  
\begin{enumerate}                                                                          
\item $(\tilde{x}_a, \lambda_{1a}, \lambda_{2a}) = (\frac{1}{2}\alpha, 1-\alpha, 0)$                                                
\item  $(\tilde{x}_b, \lambda_{1b}, \lambda_{2b}) = (1-\frac{1}{2}\alpha, 0, 1-\alpha)$                                               
\item  $(\tilde{x}_c, \lambda_{1c}, \lambda_{2c}) = (\frac{1}{2},0,0).$                          
\end{enumerate}                             

We wish to  evaluate $\dot{f}$ at all three triples to find the minimizer of $\dot{f}$. Before proceeding, we note that 
the first two triples in $\tilde{\mathcal{K}}^-_{\dot{f}}$ are defined in terms of $\alpha$, so we must evaluate feasibility for all
$\alpha \in [0,1].$  
\begin{enumerate}
\item For $(\tilde{x}_a, \lambda_{1a}, \lambda_{2a})$, all $\alpha \in [0,1]$ are feasible, and $\dot{f}(\tilde{x}_b) = \frac{1}{2}\alpha -      
\frac{1}{4}\alpha^2, \forall \alpha$.                       
\item For $(\tilde{x}_b, \lambda_{1b}, \lambda_{2b})$, all $\alpha \in [0,1]$ are feasible, and $\dot{f}(\tilde{x}_b)= \frac{1}{2}\alpha - 
\frac{1}{4}\alpha^2, \forall \alpha$.      
\item In $(\tilde{x}_c, \lambda_{1c}, \lambda_{2c})$, nothing depends on 
$\alpha$ and $\dot{f}(\frac{1}{2})=\frac{1}{4},  \forall \alpha$.                                     
\end{enumerate}

Now $(\tilde{x}_a, \lambda_{1a}, \lambda_{2a}) =  (\tilde{x}_b, \lambda_{1b}, \lambda_{2b}) = (1-\frac{1}{2}\alpha, 0, 1-\alpha)$  
is the minimizer of $\dot{f}$, so the lower endpoint of the
fuzzy interval representing $\tilde{x}(1-\tilde{x}),$ is the gradual number $z^- =  \frac{1}{2}\alpha -      
\frac{1}{4}\alpha^2.$

We repeat the steps, minimizing $-\dot{f}(\tilde{x})$, to find $z^+$, and find that $\tilde{\mathcal{K}}^+_{\dot{f}}$ consists of three triples:

\begin{enumerate}                                                                                                                              
\item $(\tilde{x}_a, \lambda_{1a}, \lambda_{2a}) = (\frac{1}{2}\alpha, \alpha-1, 0)$                                                           
\item  $(\tilde{x}_b, \lambda_{1b}, \lambda_{2b}) = (1-\frac{1}{2}\alpha, 0, \alpha-1)$                                                         
\item  $(\tilde{x}_c, \lambda_{1c}, \lambda_{2c}) = (\frac{1}{2},0,0).$                                                                        
\end{enumerate}                                                                
We wish to  evaluate $\dot{f}$ at all three triples to find the minimizer of 
$\dot{f}$. Before proceeding, we note that                       
the first two triples in $\tilde{\mathcal{K}}^-_{\dot{f}}$ are defined in terms of $\alpha$, so we must evaluate feasibility for all              
$\alpha \in (0,1].$                                                                                                                            
\begin{enumerate}                                                                                                                              
\item For $(\tilde{x}_a, \lambda_{1a}, \lambda_{2a})$, only $\alpha = 1$ is feasible, and $\dot{f}(\tilde{x}_a)= \frac{1}{4}$ for 
$\alpha = 1$ and is undefined otherwise.                                                                                                               
\item For $(\tilde{x}_b, \lambda_{1b}, \lambda_{2b})$,  only $\alpha = 1$ is feasible, and $\dot{f}(\tilde{x}_b)= \frac{1}{4}$ for               
$\alpha = 1$ and is undefined otherwise.                                                                                                                             
\item In $(\tilde{x}_c, \lambda_{1c}, \lambda_{2c})$, nothing depends on $\alpha$ and $\dot{f}(\frac{1}{2})=\frac{1}{4},~ \forall \alpha$.        \
                                                                               \end{enumerate}
                                                                               Now clearly $(\tilde{x}_c, \lambda_{1c}, \lambda_{2c}) =  (\frac{1}{2},0,0)$                
is the maximizer of $\dot{f}$, so the upper endpoint of the                                                                                    
fuzzy interval representing $\tilde{x}(1-\tilde{x}),$ is the gradual number $z^- =                                        
\frac{1}{4},~ \forall \alpha \in [0,1].$

\begin{figure}
\centerline{\includegraphics[width=3.0in]{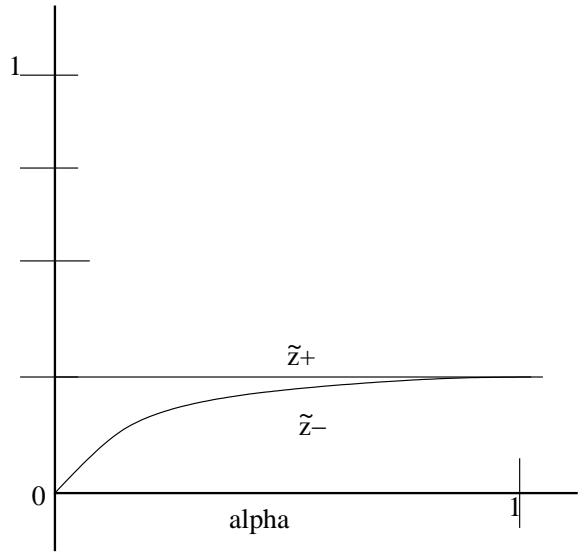}}
\caption{ $[\tilde{z}^-,\tilde{z}^+]=[\frac{1}{2}\alpha-
\frac{1}{4}\alpha^2,\frac{1}{4}],                                                 
\forall \alpha \in [0,1]$}\label{fig:fuzzyz}
\end{figure}

Our result, the fuzzy interval $[\tilde{z}^-,\tilde{z}^+]=[\frac{1}{2}\alpha-\frac{1}{4}\alpha^2,\frac{1}{4}], 
\forall \alpha \in (0,1]$ is illustrated in Figure \ref{fig:fuzzyz}.

\section{Conclusion}
Viewing a fuzzy interval as an interval of gradual numbers puts at our 
disposal a number of significant tools previously used only for real interval
analysis or for analysis of $\alpha$-cuts of fuzzy intervals.  We have shown
that optimization problems with bound constraints can be useful in 
evaluating non-monotonic functions on fuzzy intervals. 

\bibliographystyle{plain}
\bibliography{../../1thesis/latex/thesis}
\end{document}